\documentclass[a4paper,11pt]{article}
\usepackage{amssymb}
\usepackage{bm}
\usepackage{amsmath}
\usepackage{graphicx}
\usepackage{xcolor}
\usepackage{dsfont}
\usepackage{caption}
\usepackage[utf8]{inputenc}

\def\vsp{\vspace{0.2em}}
\def\vspp{\vspace{0.5em}}
\def\vsppp{\vspace{0.8em}}

\def\proof{\noindent {\bf Proof \ }}
\def\eproof{{ }\hfill$\Box$}
\def\noi{\noindent}

\def\plu{[0,\infty\mathclose[}

\def\gh{{\rm g_{\rm hyp} }}

\def\hz{{ h_z}}
\def\hzb{{ h_{\bar z}}}

\def\mk{{\cal M}_k}
\def\fk{{\cal F}_k}
\def\RR{{\mathbb R} }
\def\CC{{\mathbb C}}
\def\NN{{\mathbb N}}
\def\ZZ{{\mathbb Z}}
\def\DD{{\mathbb D} }
\def\SS{{\mathbb S} }
\newtheorem{theo}{Theorem} [section]
\newtheorem{prop}[theo]{Proposition}
\newtheorem{lemme}[theo]{Lemma}
\newtheorem{cor}[theo]{Corollary}
\newtheorem{ass}[theo]{Assertion}
\newtheorem{defi}[theo]{Definition}
\numberwithin{equation}{section}

\def\bt{\begin{theo}}
\def\et{\end{theo}}

\def\bc{\begin{cor}}
\def\ec{\end{cor}}

\def\bass{\begin{ass}}
	\def\eass{\end{ass}}

\def\bd{\begin{defi}}
\def\ed{\end{defi}}

\def\bl{\begin{lemme}}
\def\el{\end{lemme}}

\def\bp{\begin{prop}}
\def\ep{\end{prop}}

\begin{document}

\title{Harmonic quasi-isometries of pinched\\ 
Hadamard surfaces are injective}

%\title{Harmonic quasi-isometric maps IV~:\\  pinched Hadamard surfaces}

\author{Yves Benoist \&
Dominique Hulin}
\date{}

%\vfill\eject
%\setcounter{page}{1}
\maketitle

\begin{abstract} We prove that a harmonic  quasi-isometric map between pinched Hadamard surfaces is a quasi-conformal diffeomorphism. 
\end{abstract}

\renewcommand{\thefootnote}{\fnsymbol{footnote}} 
\renewcommand{\thefootnote}{\arabic{footnote}} 

%{\footnotesize \tableofcontents}

%10
\section{Introduction}\label{sec intro}

%11
\subsection{Main result}
The main result of this paper is the following. 

\bt \label{th main} Let $h:S_1\to S_2$ be a harmonic  quasi-isometric  map between pinched Hadamard surfaces. Then, $h$ is a quasi-conformal diffeomorphism.
\et

A {\it pinched Hadamard manifold} is a complete simply-connected Riemannian manifold  whose curvature satisfies $-b^2\leq K\leq -a^2$ for some positive constants $0<a\leq b$.
For instance, the hyperbolic disk $\DD$ is a pinched Hadamard surface with constant curvature $-1$. 

\vsp A map $f:M_1\to M_2$ between two metric spaces is {\it quasi-isometric} if there exists a constant $c\geq 1$ such that,  for every $x,x'\in M_1$, 
\begin{align}\label{def qi}
	c^{-1\, }d(x,x')-c\leq d(f(x),f(x'))\leq c\, d(x,x')+c.
\end{align}

\vsp
A smooth map $h:M_1\to M_2$ between Riemannian manifolds is {\it harmonic} if it is a critical point for the Dirichlet energy integral $E(h)=\int |Dh|^2 dv_{M_1}$ with respect to variations with compact support. 

\vsp 
A diffeomorphism $h:M_1\to M_2$ between $n$-dimensional Riemannian manifolds is 
{\it quasi-conformal},
if there exists a constant $C>0$ such that $\|Dh\|^n \leq C \,|{\rm Jac}(h)|$ where ${\rm Jac}(h):={\rm det}(Dh)$ is 
the Jacobian of $h$.

%12
\subsection{A few comments}
\hspace*{1em}
The special case of Theorem \ref{th main} 
where both $S_1$ and $S_2$ are the hyperbolic disk $\DD$,  
is due to Li-Tam \cite{LT} and Markovic \cite{Marko2}.

\vsppp
The main issue in Theorem \ref{th main} is the injectivity of $h$.
The quasi-confor\-mality of $h$ is but our way to prove injectivity.

\vspp
In Theorem \ref{th main} we only deal with surfaces. Indeed the analog  in higher dimension is not true. 
A counterexample due to Farrell, Ontaneda and Raghunathan is given in \cite{FOR}.

\vspp
Given two pinched Hadamard surfaces $S_1$ and $S_2$, 
there exist many harmonic quasi-isometric maps from $S_1$ to $S_2$
(see \cite{BH18} or Theorem \ref{th rappel bh18} below). 
Theorem \ref{th main} asserts that all these  maps are injective. 

\vspp 
Theorem  \ref{th main} extends the Schoen-Yau injectivity theorem in \cite{SchoenYau78}
which says that a harmonic map between two compact Riemannian surfaces 
with negative curvature, when homotopic to a diffeomorphism, is also a diffeomorphism. 
This injectivity theorem is used in the parametrization due to J. Sampson and M. Wolf of the Teichmuller space by the Hopf quadratic differentials, see \cite{Wolf89} and \cite{Jost06}.
\vspp

From a historical point of view, 
the first injectivity theorem for harmonic maps is due to 
Rado-Kneser-Choquet, almost 100 years ago. It states that, in the Euclidean plane,
the harmonic extension of an homeomorphism of the unit circle is a diffeomorphism of the unit disk,
see \cite[Lemma 5.1.10]{Hubb}. 
The analog statement in dimension $d\geq 3$ is not true. A counterexample is given by R. Laugesen in \cite{Laugesen}.
Later on, injective harmonic maps between surfaces were studied by H. Lewy in \cite{Lewy}
who proved that their Jacobian does not vanish, by R. Heinz in \cite{Heinz} and by J. Jost and H. Karcher 
in \cite[Chapter 7]{JostLNM} 
who found a lower bound for their Jacobian.
There is also an extension of the Schoen-Yau injectivity theorem
by J. Jost and R. Schoen that allows some positive curvature in \cite[Chapter 11]{JostLNM}.

%13
\subsection{Structure of the paper}
\hspace*{1em}
In Chapter \ref{sec reduction}, we recall classical facts concerning Hadamard surfaces, quasi-isometric maps and harmonic maps between surfaces.
We will see that
we can assume that the source  $S_1$ is the hyperbolic disk $\DD$. Recall that
the special case of Theorem \ref{th main} where the target $S_2$ is the hyperbolic disk $\DD$
is due to Li--Tam and Markovic.  

In Chapter \ref{sec continuity}  
we give an overview of the proof of Theorem \ref{th main}. This proof uses a deformation $(g_t)$ of the metric on $S_2$, starting with the hyperbolic metric, and a deformation $(h_t)$ of the harmonic map $h$. 
The key point will be to obtain a uniform  upper bound for the norm of the differential of $h_t$
and a uniform lower bound for the Jacobian of $h_t$.

In Chapter \ref{sec convergence}, we gather compactness results for  Hadamard surfaces and harmonic maps.

In Chapter \ref{sec quasi conf}, we obtain a uniform lower bound for the Jacobian of harmonic quasi-conformal 
diffeomorphisms.

In Chapter \ref{sec proof}, we prove that the family $(h_t)$ varies continuously with $t$ 
and we complete the proof of Theorem \ref{th main}.

In Chapter \ref{sec litam},
we include a short new proof of the special case of Theorem \ref{th main} where $S_1=S_2=\DD$. 

This paper is as self-contained as possible, the main tools being
the Bland-Kalka uniformization theorem in \cite{BlandKalka},  
the Bochner equations for harmonic maps between surfaces in \cite{Jost06},
the existence and uniqueness of quasi-isometric harmonic maps
in \cite{BH18}, and the PDE elliptic regularity in \cite{GilbargTrudinger}.

%20
\section{Background}\label{sec reduction}
\begin{quote}
We recall well-known properties of pinched Hadamard surfaces, quasi-isometric maps
and harmonic maps between surfaces. 
\end{quote}

%21
\subsection{Pinched Hadamard surfaces}
\label{sec pinched hadamard}
The first example of a pinched Hadamard surface is the hyperbolic disk  $\DD=(D , \gh ) $, where $D =\{ |z|<1\}\subset\CC$ is the unit disk equipped with the hyperbolic metric $\gh=\rho^2(z)|dz|^2$ with conformal factor  $\rho^2= 4(1-|z|^2)^{-2}$. It is a Hadamard manifold with constant curvature $-1$.

Any pinched Hadamard surface is conformal to the disk, namely reads as $(D ,\sigma^2(z) |dz|^2)$. Moreover the conformal factors $\rho^2$ and $\sigma^2$ are in a bounded ratio~:  if the curvature $K$ of this surface satisfies $-b^2\leq K\leq -a^2<0$, then $a^2\sigma^2\leq \rho^2\leq b^2\sigma^2$. See Proposition \ref{prop prescribed curvature}.

Also observe that, for maps  defined on a Riemannian surface $S_1$, 
the Dirichlet energy functional is invariant under a conformal change of metric on $S_1$. Hence,  the harmonicity of such a map depends only on the conformal class of  the source surface.

\vsp
We infer from this discussion that, to prove Theorem \ref{th main}, we 
can assume that $S_1$ is the hyperbolic disk $\DD$.

%22
\subsection{Quasi-isometric maps}
\label{sec quasi-iso}
Let $S=(D ,\sigma^2(z) |dz|^2)$ be a pinched Hadamard surface. It is a proper Gromov hyperbolic space (a general reference for Gromov hyperbolic spaces is \cite{GhysHarp90}). The boundary at infinity $\partial_\infty S$ of $S$ is defined as the set of equivalence classes of geodesic rays, where two geodesic rays are identified whenever they remain within bounded distance from each other. The union $\overline S= S\cup\partial_\infty S$ provides a compactification of $S$ (see \cite{Ball}).

\vspp
The boundary at infinity $\partial_\infty\DD$  naturally identifies with the boundary $\SS^1=\{ z\in\CC \, , \; |z|=1\}$ of $D$.  Since the identity map ${\rm Id} :D\to D$ is a quasi-isometry between the  hyperbolic disk 
$\DD=(D ,\rho^2(z) |dz|^2)$ and the surface $S=(D ,\sigma^2(z) |dz|^2)$, the boundary at infinity $\partial_\infty S$ also identifies canonically with $\partial_\infty\DD=\SS^1$.

\vspp A quasi-isometric map $f:\DD\to S$ admits a boundary value at infinity $\partial_\infty f : \partial_\infty \DD\to \partial_\infty S$, that we read as $\partial_\infty f :  \SS^1\to \SS^1$ through the above identifications. Two quasi-isometric maps share the same boundary value at infinity if and only if they remain within bounded distance from each other.
The maps $\varphi:\SS^1\to\SS^1$ that appear as boundary values at infinity of quasi-isometric maps $f:\DD\to S$ are exactly the quasi-symmetric homeomorphisms. 
For convenience, we  identify $\SS^1$ with $\RR/2\pi\ZZ$.

\bd 
\label{def quasi sym} 
Let $k\geq 1$. An homeomorphism $\varphi :\SS^1\to\SS^1$ is a $k$-quasi-symmetric map  if 
\begin{equation}\label{eq quasi sym}
	\frac 1 k\leq \frac  {\varphi (\theta+\alpha)-\varphi (\theta)}     {\varphi (\theta)-\varphi (\theta-\alpha)}\leq k
\end{equation}
holds for every $\theta , \alpha$ with $0<\alpha\leq\pi$.
\ed

Note that any quasi-isometric map $f:\DD\to S$ is actually a quasi-isometry. Namely, there exists
$C>0$ such that $d(y,f(\DD))\leq C$ holds for all $y$ in $S$.  Indeed,
the inverse $\varphi^{-1}$ of its boundary map is also a quasi-symmetric homeomorphism, 
hence $\varphi^{-1}$ is the boundary map 
of a quasi-isometric map $f':S\to\DD$, and the map $f\circ f':S\to S$ is within bounded
distance from the identity map. 

In a previous paper, we studied harmonic quasi-isometric maps between pinched Hadamard manifolds. Our result, when specialized to surfaces, asserts that any quasi-isometric map $f:\DD\to S$ has the same boundary value at infinity as a unique harmonic quasi-isometric map. In other words,  the following holds. 
\bt\label{th rappel bh18} {\rm \cite{BH18}  }Let $S=(D ,\sigma^2(z) |dz|^2)$ be a pinched Hadamard surface and $\varphi :\SS^1\to\SS^1$ be a quasi-symmetric map. Then, there exists a unique harmonic quasi-isometric map $h:\DD\to S$ such that $\partial_\infty h =\varphi$.
\et

%23
\subsection{Harmonic maps between surfaces}
\label{sec harmonic}

\begin{quote}
We introduce some notation that will be used throughout the paper, and recall some classical results concerning harmonic maps between surfaces.  A general reference for this section is Jost \cite{Jost06}.
\end{quote}
Let $h:\DD\to S$ be a smooth map from the hyperbolic disk  $\DD=(D ,\rho^2(z)|dz|^2)$  to a pinched Hadamard surface $S=(D ,\sigma^2(z) |dz|^2)$ with pinching condition $-b^2\leq K\leq -a^2<0$.
Recall that the curvature $K$ of $S$ is given by 
$$
K=-\sigma^{-2}\,\Delta_e\log\sigma
$$
where $\Delta_e=4\partial_z\partial_{\overline{z}}$ is the Euclidean Laplacian.
For such a map $h$, we introduce as usual the functions $\hz,\hzb:\DD\to\CC$ defined by
$$\hz=\frac 12(h_x-ih_y), \quad \hzb=\frac 12(h_x+ih_y)\, $$
where the conformal parameter reads as $z=x+iy$, and the subscript $x$ or $y$ indicates a directional derivative. The map $h$ is holomorphic (or anti-holomorphic) if $\hzb=0$ (or $\hz=0$). It is worth noting that $\overline{\hzb}=\bar h_z$.

\bp\label{prop def harmonic}
{\rm \cite[Section 3.6]{Jost06}}
The map $h:\DD\to S$ is harmonic if and only if it satisfies
\begin{equation*}
	h_{z\bar z}+2\,(\frac {\sigma_z}{\sigma}\circ h) \,\hz\hzb=0\, .
\end{equation*}
\ep
If the map $h$ is either holomorphic, or anti-holomorphic, then it is harmonic.
Introduce the square norms of the complex derivatives of $h$~:
\begin{eqnarray*}
	H=\| \partial h\|^2:=\frac{\sigma^2\circ h}{\rho^2}\, |\hz|^2\quad \hbox{and}\quad
	L=\|\bar\partial h\|^2:=\frac{\sigma^2\circ h}{\rho^2}\, |\hzb|^2\, ,
\end{eqnarray*}
so that one has $\|Dh\|^2=H+L$.
Observe that $h$ is a local diffeomorphism if the Jacobian $J=H-L$ does not vanish, and is moreover orientation preserving if $J>0$.

\bl
\label{lem bochner} 
{\rm\cite[Section 3.10]{Jost06}}  Let $h:\DD\to S$ be a harmonic map. On the open subsets where they are non zero, the functions $H$ and $L$ satisfy the Bochner equations
\begin{eqnarray}
\label{eq logh}
(1/2)\, \Delta\log H&=& \!(-K\circ h)\, J-1 \, , \\
\label{eq logl}
(1/2)\, \Delta\log L&=& \; (\, K\circ h)\, J-1\, .
\end{eqnarray}
\el
Here $\Delta = 4\,\rho^{-2}\partial_z\partial_{\bar z}$ is the Laplace operator relative to the hyperbolic metric.

On the open set $\Omega:=\{ \hz\neq 0\}$, 
we introduce the conformal distortion 
$\mu :\Omega\to\CC$ by letting $\hzb =\mu\hz$,
so that one has the useful  equalities  
\begin{eqnarray}
\label{eq lsurh jsurh}
|\mu|^2=L/H
\;{\rm ,}\;\; 
1-|\mu|^2=J/H\, 
\;\;{\rm and}\;\;\; 
\frac{1-|\mu|^2}{1+|\mu|^2}=\frac{J}{\|Dh\|^2}\, .
\end{eqnarray}

%30
\section{A family of metrics and harmonic maps}
\label{sec continuity}
\begin{quote}
In this  section we explain the continuity method that will be used 
to prove Theorem \ref{th main}.
\end{quote}

Let $S=(D ,\sigma^2(z) |dz|^2)$ be a pinched Hadamard surface, with curvature bounds $-b^2\leq K\leq -a^2<0$. Choose an increasing quasi-symmetric homeomorphism $\varphi : \SS^1\to\SS^1$, and let $h:\DD\to S$ be the unique harmonic  quasi-isometric map with boundary value at infinity $\partial_\infty h=\varphi$. 
We want to prove that $h$ is a quasi-conformal diffeomorphism. 

\vspp In case the surface $S$ is the hyperbolic disk, that is for a harmonic quasi-isometric map $h:\DD\to\DD$, the result is due to Li-Tam and Markovic (see Chapter \ref{sec litam} for a proof). 
To prove it for a harmonic map $h:\DD\to S$ with values in a general pinched Hadamard surface $S$, we use the method of continuity, involving a  family of pinched Hadamard surfaces $S_t=(D , e^{2u_t}\gh)$, for $0\leq t\leq 1$, starting with $S_0=\DD$ and such that $S_1=S$.

%31
\subsection{Construction of the metrics $g_t$}
\begin{quote}
We construct the metric $g_t$ by prescribing its curvature.
\end{quote}

More specifically, we introduce  for $0\leq t\leq 1$  
the unique complete  conformal metric $g_t=e^{2u_t}\gh$ 
on the unit disk $D$ with curvature 
$K_t:=\! -(1\! -\! t)\!+\! tK$. Each function $K_t$ being pinched between 
two negative constants, the exis\-tence and uniqueness of such a metric 
is granted by the following.

\bp
\label{prop prescribed curvature} 
{\rm\cite{BlandKalka} }
Let $k$ be a smooth function on the unit disk $D$ such that 
$-\beta^2\leq k\leq -\alpha^2$ for some constants $0<\alpha\leq \beta$. 
Then, there exists a unique complete 
conformal metric $ g=e^{2u}\gh$ on $D$ with curvature $k$. 
Moreover, the conformal factor $e^{2u}$ 
is controlled, with $\beta^{-2}\leq e^{2u}\leq \alpha^{-2}$.
\ep

\noi
We do not reproduce here the proof that is given in \cite{BlandKalka} and that relies on the sub-supersolution method for the curvature equation 
\begin{equation}\label{eq curvature}
	\Delta  u = (- k)\, e^{2u} -1\, ,
\end{equation}
where, as above,
$\Delta$ is the Laplace operator for the hyperbolic metric $\gh$.

The proof also uses  the generalized maximum principle of Yau in  \cite{Yau75}.
We will need  later a light form of this principle that reads as follows.

\bl
\label{lem max yau}
Let $v:S\to\RR$ be a smooth function defined on a pinched Hada\-mard surface $S$. 
Assume that $v$ is bounded above. 

Then, there exists a sequence $(x_n)$ in $S$  such that
\begin{eqnarray}
\label{eq yau}
\mbox{$v(x_n)\to\sup _Sv$,\;\; $|\nabla v|(x_n)\to 0$ and\;\; 
$\displaystyle\limsup\Delta v(x_n)\leq 0$.}
\end{eqnarray}
\el

\proof
We can assume that $\sup_Sv=1$. 
We fix a point $x_0\in S$ where this supremum is not achieved and we introduce
the function $v_n$ on $S$ given by $v_n(x)=v(x)\,e^{-d(x,x_0)/n}$. 
This function is smooth, except maybe at $x_0$, and
it achieves its supremum  at a point $x_n\neq x_0$ for $n$ large. 
This sequence $(x_n)$ satisfies \eqref{eq yau} since
$v_n(x_n)\to 1$,\; $\nabla v_n(x_n)=0$ and \; $\Delta v_n(x_n)\leq 0$.
\eproof

%32
\subsection{Construction of the harmonic maps $h_t$}
\hspace*{1em}
We construct the harmonic map $h_t$ by prescribing its boundary map.
\vspp 

By construction, one has $\DD=(D,g_0)$ and $S=(D,g_1)$. 
For $0\leq t\leq 1$, we let $h_t:\DD\to S_t$ be the unique harmonic quasi-isometric map 
whose boundary value at infinity is $\varphi :\SS^1\to\SS^1$. 
Recall that the existence and uniqueness of those $h_t$ 
are granted by Theorem \ref{th rappel bh18}.

\vspp
Here are some basic information concerning these harmonic maps $h_t$. 
For $0\leq s,t\leq 1$, let 
$\displaystyle d(h_s,h_t):=\sup_{z\in D} d(h_s(z),h_t(z))$ denote 
the uniform distance between these two maps, where the distance is taken 
with respect to the hyperbolic metric $g_{\rm hyp}$ on the target. 

\bl
\label{lem hs ht}
There exists $c_*>0$ such that, for all $t\in [0,1]$,
the map $h_t$ is $c_*$-quasi-isometric, one has $d(h_t,h_0)\leq c_*$,
and the map $h_t$ is $c_*$-Lipschitz. 
\el

\noi
Remark that, since the functions  $u_t$ are uniformly bounded (Proposition \ref{prop prescribed curvature}), it was not really necessary to specify with respect to which one of the metrics $g_t$ the above
distances were being estimated.

\vspp
\proof 
As explained in Section \ref{sec quasi-iso}, there exists 
a $c$-quasi-isometric map $f:\DD\to\DD$
whose boundary value at infinity is our quasi-symmetric map $\partial_\infty f=\varphi$.  
By taking a larger  constant $c$, we may assume  that each map $f:\DD\to S_t$ (that is, the same map $f$ now seen with values in one of the Riemannian surfaces $S_t$, with $t\in [0,1]$) is  $c$-quasi-isometric. 

Thus the main result of \cite{BH18}  asserts that there exists a constant $C>0$  such that $d(f,h_t)\leq C$.
This constant $C$ depends only on $c$ and on the pinching constants $a$ and $b$, 
hence it does not depend on $t\in [0,1]$. Thus the first two claims hold if $c_*\geq 2c+2C$.

The map $f$ being $c$-quasi-isometric, each harmonic map $h_t:\DD\to S_t$ sends any ball $B(z,1)\subset\DD$ with radius $1$ inside the ball $B(h_t(z),R)\subset S_t$ with radius $R=2c+2C$. Now the 
uniform Lipschitz continuity of the maps $h_t$ follows from the Cheng lemma, that we recall below. \eproof

\bl 
{\rm\cite{Cheng80}}
\label{lem cheng} 
Let $S$ be a  Hadamard surface with $-b^2\leq K\leq 0$. 
There exists a constant $\kappa$, that depends only on  $b$,  
such that if a harmonic map $h:\DD\to S$ satisfies $h(B(z,1))\subset B(h(z),R)$ for some radius $R$, then
$$
\| Dh(z)\|\leq\kappa\, R\, .
$$
\el

%33
\subsection{An injectivity criterion}
The following lemma tells us that a uniform lower bound for the Jacobian $J_t={\rm Jac}(h_t)$
is enough to ensure that 
$h_t$ is a quasi-conformal diffeomorphism.

\bl
\label{lem continuity methodB}  
If  $ \displaystyle \inf_{z\in \DD}J_t(z)>0$ then $h_t$ is a quasi-conformal diffeomorphism.
\el

\proof 
By assumption the Jacobian $J_t$ does not vanish, hence the map $h_t:\DD\to S_t$ is a local diffeomorphism. By  construction, the map $h_t$ is quasi-isometric, hence it is a proper map. It thus follows that
$h_t$ is a covering map. Hence, since $S$ is simply connected, the map $h_t$ is a diffeomorphism. 
Since, by Lemma \ref{lem cheng}, $h_t$ is Lipschitz, the lower bound for its Jacobian $J_t$ ensures that $h_t$ is quasi-conformal.
\eproof

%34
\subsection{Strategy of proof of Theorem \ref{th main}}
\label{sec strategy} 

We will need  the following two propositions.

\bp
\label{prop continuity methodA}  
There exists  $j_*>0$ such that, for all $t\in [0,1]$ for which $h_t$ is a quasi-conformal 
diffeomorphism, one has  $J_t\geq j_*$.
\ep

Proposition \ref{prop continuity methodA} 
is a straightforward consequence of  Proposition \ref{prop delta de m}
that will be proven in Chapter \ref{sec quasi conf}. 
Indeed Lemma \ref{lem hs ht} ensures that the 
maps $h_t$ are $c_*$-Lipschitz.

Let $C_b(\DD,\RR)$ be the space of bounded continuous functions $\psi$  endowed with the sup norm~: 
$\|\psi\|_\infty= \displaystyle\sup_{z\in\DD}|\psi(z)|$.

\bp
\label{prop continuity methodB} 
The map $t\in [0,1]\to J_t \in C_b(\DD,\RR)$ is continuous. 
\ep

Proposition \ref{prop continuity methodB} will be proven in Chapter \ref{sec proof}
as part of Proposition \ref{prop continuity methodBplus} . 
\vspp 

\noi{\bf Proof of Theorem \ref{th main} using Propositions \ref{prop continuity methodA} and
\ref{prop continuity methodB} }
Let $A$ be the set of  parameters $t\in [0,1]$ such that 
the harmonic map $h_t:\DD\to S_t$ is a quasi-conformal diffeomorphism. 
We want to prove that $1\in A$.
We already know that $0\in A$ (this is Theorem \ref{th litam} due to Li-Tam and Markovic).
It is enough to check that $A$ is open and closed.
Let $j$ be the function  on $[0,1]$ given by 
$$
j(t):=\displaystyle \inf_{z\in \DD}J_t(z)\in \RR.
$$
By Proposition \ref{prop continuity methodB} the function 
$j$ is continuous.
By Lemma \ref{lem continuity methodB} and Proposition \ref{prop continuity methodA}, one has both $A=j^{-1}(]0,\infty[)$ and $A=j^{-1}([j_*,\infty[)$. Hence $A$ is both open and closed. \hfill$\Box$

%40
\section{Sequences of metrics and harmonic maps}\label{sec convergence}
\begin{quote}
In order to obtain the uniform lower bounds in Chapter \ref{sec quasi conf}, 
or the continuity properties in Chapter \ref{sec proof}, 
we will  have to consider sequences of conformal metrics on the unit disk $D$, 
and sequences of harmonic maps. 
In this chapter, we state compactness results for such sequences.
\end{quote}

\noi
These compactness results also hold in higher dimension (see \cite{Petersen16}, or \cite{BH18}). 
Since we will only deal here with conformal metrics on the disk $D$, the complex parameter $z\in D$ naturally provides a global harmonic chart for these metrics 
so that the statements and the proofs are more elementary.

%41
\subsection{Sequence of Hadamard surfaces}
\begin{quote}
Let us begin with sequences of conformal Riemannian structures on the unit disk $D$. 
\end{quote}
Convergence in the following lemma is a special case of the Gromov-Hausdorff convergence 
for isometry classes of pointed proper metric spaces using the base point $0\in D$. 
See   \cite[\S 5.3]{BH18} or \cite{BuragoBuragoIvanov} for a short introduction to this notion. 

\bl
\label{lem seq surfaces} 
Let  $g_n=e^{2u_n}\gh$ 
be a sequence of complete conformal  metrics  on the unit disk $D$ with curvature $-b^2\leq K_n\leq -a^2<0$. 
Then there is a subsequence of 
$(u_n)$ that converges to a $C^1$ function $u_\infty$ in the $C^1_{\rm loc}$ topology. 

The limit metric $g_\infty=e^{2u_\infty}\gh$ is a $C^1$ complete conformal metric on  $D$,
and $S_\infty:=(D,g_\infty)$ is a CAT-space with curvature between  $-b^2$ and $-a^2$.
\el

\proof  
Proposition \ref{prop prescribed curvature} ensures that 
the logarithms $u_n:D\to\RR$  of the conformal factors are uniformly bounded.  
The curvature equation 
\begin{equation}
	\Delta  u_n = (- K_n)\, e^{2u_n} -1\,  \tag{\ref{eq curvature}${}_n$}
\end{equation}
for $g_n$ ensures that the Laplacians $\Delta u_n$ are also uniformly bounded. 

Pick $0 \leq \alpha <1$.  We may now apply to the sequence $(u_n)$ the following first Schauder estimates (see \cite[Theorem 3.9]{GilbargTrudinger} or \cite[Theorem 70]{Petersen16}). 
These estimates state that there exists a constant $c$ such that, 
for any smooth function $v:\DD\to\RR$ on the hyperbolic disk, the inequality
\begin{equation}\label{eq schauder 1}
	\| v\|_{C^{1,\alpha}(B_1)}\leq c_\alpha (  \| \Delta v\|_{C^{0}(B_2)}+\| v\|_{C^{0}(B_2)})\, 
\end{equation}
holds for any pair of concentric hyperbolic balls $B_1\subset B_2\subset\DD$ with respective  radii $1$ and $2$.
This provides a uniform local bound for the norms $\| u_n\|_{C^{1,\alpha}}$. Going if necessary to a subsequence,  we may thus assume that the sequence $(u_n)$ converges in the $C^1_{\rm loc}$ topology. Let  $u_\infty=\lim u_n$ and  $g_\infty=e^{2u_\infty}\gh$ and introduce $S_\infty =(D,g_\infty)$.
As a limit of such, the length space $S_\infty$ is a CAT-space with curvature between  $-b^2$ and $-a^2$ (see  \cite[Corollary II.3.10]{BridsonHaefliger} 	and \cite[Theorem 10.7.1]{BuragoBuragoIvanov}).
\eproof
\vspp 

\noi{\bf Remark } Under the hypothesis of Lemma \ref{lem seq surfaces},
after extraction, the sequence of bounded functions $K_n:D\to\RR$ converges   weakly
to a bounded measurable function $K_\infty:D\to\RR$ with $-b^2\leq K_\infty\leq -a^2$,
and the $C^1$ function $u_\infty$ is a weak solution of
\begin{equation}
\Delta  u_\infty = (- K_\infty)\, e^{2u_\infty} -1\, . \tag{\ref{eq curvature}${}_\infty$}
\end{equation}

%42
\subsection{Sequence of harmonic maps}
\hspace*{1em}
Now turn to sequences of maps between such Riemannian surfaces.
\bl
\label{lem seq maps}  
Let $S_n=(D,g_n)$ be a sequence converging to $S_\infty=(D,g_\infty)$ as in Lemma \ref{lem seq surfaces}. Let $c>0$, 
and let $h_n:  \DD\to  S_n$ be $c$-Lipschitz maps satisfying $d_n(h_n(0),0)\leq c$. 
Then there is a subsequence of $(h_n)$ that converges locally uniformly to a $c$-Lipschitz map
$h_\infty: \DD\to  S_\infty\, .$\\
\noi
$a)$ If all the maps $h_n$ are $C$-quasi-isometric, then $h_\infty$ is $C$-quasi-isometric.\\
$b)$ If  all the maps $h_n$ are harmonic, then $h_\infty$ is $C^2$ and is harmonic too.
\el

\proof 
Observe that, on any fixed compact set, the maps $h_n : \DD\to  S_\infty$ are $c_n$-Lipschitz for some constants $c_n$ converging to $c$.
Indeed these are the initial maps $h_n$, albeit with the limit metric on the target.
Since we assumed that $d_n(h_n(0),0)\leq c$, 
these maps $h_n$ are locally uniformly bounded 
(this means locally in $z$ and uniformly in $n$).
It thus follows from the Ascoli lemma that we may assume  
the sequence $(h_n)$ to converge uniformly on compact sets
to a $c$-Lipschitz map $h_\infty:  \DD\to S_\infty$.

$a)$ If all  $h_n:  \DD\to  S_n$ are $C$-quasi-isometric,  
then, on any fixed compact set, 
the maps $h_n : \DD\to  S_\infty$ are $C_n$-quasi-isometric 
for some constant $C_n$ converging to $C$, and so $h_\infty$
is $C$-quasi-isometric.
\vsp

$b)$ Now assume  that each map $h_n:\DD\to S_n$ is harmonic, 
namely that each function $h_n:\DD\to D\subset\CC$ satisfies  the equation
\begin{equation}
\label{eq fn harmonic}
	(h_n)_{z\bar z}+2\, ((u_n)_z\circ h_n) \,(h_n)_z(h_n)_{\bar z}=0\, .
\end{equation}
We want to prove that $h_\infty$ is harmonic, namely that it is $C^2$ and satisfies
\begin{equation}
\label{eq fi harmonic}
	(h_\infty)_{z\bar z}+2\, ((u_\infty)_z\circ h_\infty) \,(h_\infty)_z(h_\infty)_{\bar z}=0\, .
\end{equation}
The maps $h_n:\DD\to S_n$ are $c$-Lipschitz, so that all the derivatives $(h_n)_z$ and $(h_n)_{\bar z}$ are locally uniformly bounded. 
We have seen in the proof of Lemma \ref{lem seq surfaces} that the gradients $\nabla u_n$
are locally uniformly bounded, hence $(u_n)_z\circ h_n$ are  locally uniformly bounded.
Then \eqref{eq fn harmonic} ensures that the functions $\Delta h_n$ are also
locally uniformly bounded. 
We apply the first Schauder  estimates \eqref{eq schauder 1} to the functions $v=h_n$.
This implies that, for $0<\alpha<1$, 
the functions $h_n$ are uniformly bounded in the $C^{1,\alpha}_{\rm loc}$ topology.

Plugging this information in \eqref{eq fn harmonic}, 
and remembering from the proof of Lemma \ref{lem seq surfaces} that the gradients $\nabla u_n$
are also uniformly bounded in the $C^\alpha_{\rm loc}$ topology, we see that the 
functions $\Delta h_n$ are uniformly bounded in the $C^{\alpha}_{\rm loc}$ topology.
We will now apply  the second Schauder estimates to the functions $v=h_n$
(see \cite[Theorem 70]{Petersen16}). 
With the same notation as \eqref{eq schauder 1},
these estimates state
\begin{equation}\label{eq schauder 2}
\| v\|_{C^{2,\alpha}(B_1)}\leq 
c_\alpha \bigl(  \| \Delta v\|_{C^{\alpha}(B_2)}+\| v\|_{C^{0}(B_2)}\bigr).
\end{equation}
Hence the functions $h_n$ are uniformly bounded in the $C^{2,\alpha}_{\rm loc}$ topology. 

Therefore $(h_n)$ admits a  subsequence which converges in the $C^2_{\rm loc}$ topo\-logy. 
This proves that $h_\infty$ is $C^2$ and going to the limit in \eqref{eq fn harmonic} ensures that the limit map $h_\infty$ is harmonic, as claimed.
\eproof

%50
\section{A lower bound for the Jacobian}
\label{sec quasi conf}
\begin{quote}
In this section we provide a lower bound for the Jacobian $J_t$ of $h_t$ when 
$h_t$ is a quasi-conformal 
diffeomorphism (Proposition \ref{prop continuity methodA}).
\end{quote}

\noi
The notation are those of Section \ref{sec harmonic}~: 
$S$ is a pinched Hadamard surface and $h:\DD\to S$
is an harmonic map.
We assume moreover that $h$ is an orientation preserving  diffeomorphism. 
The  Jacobian of $h$, which is $J=H-L$ with  $H:=\|\partial h\|^2$ and $L:=\|\overline{\partial} h\|^2$,  
is positive.
The function $w:=\frac12\log H$ satisfies Equation \eqref{eq logh}, that we may also write as
\begin{equation}\label{eq deltaw}
	\Delta w=	(-K\circ h) \, (1-|\mu|^2)\, e^{2w}-1\, ,
\end{equation}
where $\mu:=\hzb/\hz$ is the conformal distortion.
By \eqref{eq lsurh jsurh} the diffeomorphism $h$ is quasi-conformal if and only if there exists a  $\delta<1$ such that  $|\mu|\leq\delta$.

%51
\subsection{Controlling the norm of the differential}
The next lemma tells us that the norm of the differential  $\| Dh\|$ 
of a harmonic quasi-conformal diffeomorphism is uniformly bounded below (see  also \cite{Wan92}).

\bl
\label{lem H fonction delta} 
Let $h:\DD\to S$ be a  quasi-conformal harmonic diffeomorphism, where $S$  is a pinched Hadamard surface with curvature  $-b^2\leq K\leq -a^2<0$. 
Then one has $e^{2w}\geq b^{-2}\, .$
\el

\proof
Introduce the conformal metric $\tilde g=e^{2w}\gh$ on $D$. We first prove that $\tilde g$ is complete with pinched negative curvature.  Proposition \ref{prop prescribed curvature} will then provide the lower bound on $w$.

Let $S=(D,\sigma^2(z)|dz|^2)$. The map $h:\DD\to S$ being a diffeomorphism and $S$ being complete, the pull back metric $G=h^*(\sigma^2(z) |dz|^2)$ is complete. 
This pull-back metric reads as
$G=(\sigma^2\circ h)|\hz|^2|dz\! +\!\mu d\bar z|^2$.
Since one has $\tilde g=(\sigma^2\circ h)|\hz|^2|dz|^2$ and $|\mu|\leq 1$,  one easily checks  that 
$G\leq 4\tilde g$. This ensures that the metric $\tilde g$ is complete. 

Comparison of Equation \eqref{eq deltaw} satisfied by $w$ and the curvature equation \eqref{eq curvature} yields that the  metric $\tilde g$ has curvature $\tilde K=(K\circ h)(1-|\mu|^2)$.
It follows that $-b^2\leq \tilde K\leq -a^2(1-\delta^2)<0$,
where $\delta:=\|\mu\|_\infty <1$. Proposition \ref{prop prescribed curvature} 
thus ensures that $w$ satisfies
 $b^{-2}\leq e^{2w}\leq a^{-2}(1-\delta^2)^{-1}\, .$ 
\eproof

%52
\subsection{Controlling the Jacobian}

The following proposition tells us that the Jacobian 
of a harmonic quasi-conformal diffeomorphism is controlled by its Lipschitz constant.

\bp
\label{prop delta de m} 
Let $0<a\leq b$.
Then, for every  $c>0$, there exists $j_*=j_*(a,b,c)>0$  such that, if $S$ is a pinched Hadamard surface with curvature $-b^2\leq K\leq -a^2$,  the Jacobian $J$ of any $c$-Lipschitz quasi-conformal harmonic diffeomorphism $h:\DD\to S$ satisfies $J\geq j_*$.
\ep

\proof  Assume by contradiction that  there exist a sequence of pinched Hadamard surfaces
$S_n=(D,e^{2u_n}\gh)$
with curvatures $-b^2\leq  K_n\leq -a^2$,  
a sequence $h_n:\DD\to S_n$  of $c$-Lipschitz harmonic quasi-conformal diffeomorphisms 
and  a sequence $(x_n)$ of points of $D$ such that
the Jacobian $J_n$ of $h_n$ satisfy $J_n(x_n)\to 0$.

Choosing sequences $(\gamma_n)$ and $(\gamma'_n)$ of isometries of the hyperbolic disk
such that $\gamma_n(x_n)=0$ and $\gamma'_n(h_n(x_n))=0$, and
replacing $u_n$ by $u_n\circ {\gamma'_n}^{-1}$ and  $h_n$ by $\gamma'_nh_n\gamma_n^{-1}$,
we can assume that $x_n=0$ and $h_n(x_n)=0$.

By Lemmas \ref{lem seq surfaces} and \ref{lem seq maps}, going to a subsequence, one may assume that~:\\
-- the sequence $(u_n)$ converges  to a $C^1$ function $u_\infty $
in the $C^1_{\rm loc}$ topology.\\
-- the sequence $(h_n)$ converges  to a $C^2$ map $h_\infty$ in the $C^2_{\rm loc}$ topology.

Recall from \eqref{eq lsurh jsurh} that $J_n=(1-|\mu_n|^2)e^{2w_n}$
where $\mu_n=(h_n)_{\bar z}/(h_n)_z$ is the conformal distorsion and 
where $e^{2w_n}=\|\partial h_n\|^2$. 
Lemma \ref{lem H fonction delta} ensures that
\begin{equation}\label{eq 2w}
e^{2w_\infty}=\lim_{n\to\infty}e^{2w_n}\geq b^{-2}\, .
\end{equation}
Thus $(h_\infty)_z$ does not vanish.
Hence the functions  $\mu_n$ 
also converge to a $C^1$ functions  $\mu_\infty$ in the $C^1_{\rm loc}$ topology,
and one has $\|\mu_\infty\|_\infty=1$ and $|\mu_\infty(0)|=1$.
\vspp

{\bf First step } We claim that $|\mu_\infty|\equiv 1$.

\noi
Indeed, we introduce the non negative $C^1$ functions $\ell_n:=-\log |\mu_n|^2$ defined on $\Omega_n:=\{ \mu_n\neq 0\}$ and their limit $\ell_\infty:=-\log |\mu_\infty|^2$, 
which is defined on $\Omega_\infty:=\{ \mu_\infty\neq 0\}$. 
By assumption, the function $\ell_\infty$ is a non-negative function that achieves its minimum $\ell_\infty(0)=0$ 
at the origin. 
We will prove that the set $\{ \ell_\infty =0\}$ is open in $\Omega_\infty$, so that 
$\ell_\infty\equiv 0$ as claimed. 

\vsp 
The function $\ell_n$ satisfies the equation on $\Omega_n$,
difference of \eqref{eq logh} and \eqref{eq logl}:
\begin{equation}
\label{eq mu}
 \Delta \ell_n= 4\,(-K_n\circ h_n)\, (1-e^{-\ell_n})\, e^{2w_n}\, .
\end{equation} 
Since  $|K_n|\leq b^2$,\; $1\! -\! e^{-\ell_n}\leq \ell_n$ and $e^{2w_n}\leq c^{2}$,  we infer that
\begin{equation*}\label{eq elln}
	\Delta \ell_n\leq 4b^2 c^2\, \ell_n\, .
\end{equation*}
Hence
$\ell_\infty$ is a $C^1$ function on $\Omega_\infty$ that 
satisfies in the weak sense
\begin{equation*}\label{eq ell}
	\Delta \ell_\infty\leq 4b^2 c^2\, \ell_\infty\, .
\end{equation*}
In particular, one has bounds  $\Delta_e\ell_\infty\leq C_K\ell_\infty$ on compact sets $K$  of $\Omega_\infty$ and, by 
Lemma \ref{lem heinz} below, the set $\{\ell_\infty=0\}$ is open. This proves  $|\mu_\infty|\equiv 1$.
\vspp 

{\bf Second step } We reach a contradiction.

\noi
We  recall that the functions $w_n$ satisfy \eqref{eq deltaw}, namely
$$
\Delta w_n=	(-K_n\circ h_n) \, (1-|\mu_n|^2)\, e^{2w_n}-1\, . 
$$
Since the functions $(-K_n\circ h_n)$ and $  e^{2w_n}$ are uniformly bounded and since
$\lim_{n\to\infty}|\mu_n|=1$, 
the limit function  $w_\infty=\lim w_n$ satisfies 
$\Delta w_\infty =-1$ in the weak sense. In particular $w_\infty$ is smooth.
Note also that \eqref{eq 2w} yields the lower bound $w_\infty\geq \log b^{-2}$.
 
In conclusion, $w_\infty$ is a smooth function on $\DD$ which is bounded below  
and satisfies $\Delta w_\infty=-1$.
By the generalized maximum principle 
of Lemma \ref{lem max yau}, such a function $w_\infty$ does not exist.
Contradiction.
\eproof
\vspp

In the previous proof we have used the following lemma as in \cite{Heinz}.  

\bl
\label{lem heinz} 
Let $C>0$ and $\ell$ be a non-negative  continuous function on 
an open set $U\subset\RR^2$ such that 
$\Delta_e\ell \leq C\ell$ weakly. Then the set $\{ \ell=0\}$ is open.
\el

\proof  We can assume that $\ell(0)=0$.
By a standard convolution argument, in a small ball $B(0,R)\subset \Omega$, we can  write $\ell$ as a uniform limit 
of non-negative $C^2$-functions $\ell_n$ that also satisfy 
\begin{eqnarray} 
\label{eq pour heinz}
\Delta_e\ell_n\leq C\ell_n\, .
\end{eqnarray} 
We introduce   the mean values of $\ell_n$ and $\ell$ on circles of radius $r\leq R$, 
$$
\hbox{$M_n(r):=\frac{1}{2\pi}\int_0^{2\pi}\ell_n (r\, e^{i\theta })\, d\theta$ \ and\  $ M(r):=\frac{1}{2\pi}\int_0^{2\pi}\ell(r\, e^{i\theta })\, d\theta$}.
$$
The Green representation formula (see H\"{o}rmander \cite[p.119]{Hormander}) gives
$$
\ell_n (0)= M_n(r)-\frac{1}{2\pi}\int_{B(0,r)}\Delta_e\ell_n (y)\, \log\frac{r}{|y|}\,  dy\, .
$$
Since $\ell_n$ converges uniformly to $\ell$ and $\ell (0)=0$ we infer, using \eqref{eq pour heinz}, that
$$
M(r)\leq \frac{C}{2\pi}\int_{B(0,r)}\ell (y)\,  \log\frac{r}{|y|}\, dy\, ,
$$
so that, for every $r\leq R$,
$$
M(r)\leq \frac{C\, R^2}{4}\; \sup_{[0,R]}M(t) \, .
$$
Choosing $R^2<4/C$, we obtain that $\ell\equiv 0$ on the ball $B(0,R)$.
\eproof

%60
\section{Continuity of the Jacobian}
\label{sec proof}
\begin{quote}
In this  section we prove that the metrics $g_t$,  the harmonic maps $h_t$
and their  Jacobians $J_t$ depend continuously on $t$, thus proving
Proposition \ref{prop continuity methodB}. 
\end{quote}

%61
\subsection{A continuous family of metric}
In Chapter \ref{sec continuity}, we introduced   pinched Hadamard surfaces $S_t=(D,e^{2u_t}\gh)$ with curvature $K_t=(t-1)+tK$, where $-b^2\leq K\leq -a^2<0$ ($t\in[0,1]$). In particular, $S_0=\DD$.
We have seen that all the metrics $g_t$  are uniformly bi-Lipschitz to each other.   This means that the functions $u_t:D\to\RR$  are uniformly bounded.

Lemma \ref{lem u borne et continu} tells us that they are uniformly bounded in norm $C^{1}$
and that the map $t\in [0,1]\to u_t\in C^1$ is continuous. 
Here the gradients $\nabla $, as well as their norms, are taken with respect to the hyperbolic metric $\gh$. 

\bl 
\label{lem u borne et continu} 
There exists a constant $c$ such that, for every $0\leq t\leq 1$
\begin{eqnarray}
\label{eq u et grad u}	\| u_t\|_\infty+\| \nabla u_t\|_\infty&\leq&c\\
\label{eq u et grad u continu}		
\| u_t-u_s\|_\infty+\| \nabla (u_t-u_s)\|_\infty&\leq& c\, |t-s|\, .
\end{eqnarray}
\el

\proof  
We argue as in the proof of Lemma \ref{lem seq surfaces}.
Let us first prove \eqref{eq u et grad u}.
Each conformal factor $e^{2u_t}$ is solution of the curvature equation \eqref{eq curvature}, here
\begin{equation}\label{eq ut}
	\Delta u_t= (-K_t)\, e^{2u_t}-1\, .
\end{equation}
Since the metrics $g_t$ are complete, and the $K_t$ satisfy a uniform pinching condition $-B^2\leq K_t\leq -A^2<0$ for all $0\leq t\leq 1$, Proposition \ref{prop prescribed curvature} ensures that the  functions $u_t$ are uniformly bounded.
Plugging into \eqref{eq ut}, we infer  that the Laplacians $\Delta u_t$  are also uniformly bounded.
Hence the Schauder estimates \eqref{eq schauder 1} with $\alpha=0$
and $v=u_t$ yield the uniform bound \eqref{eq u et grad u}.

\vspp  
We now prove \eqref{eq u et grad u continu}. 
Using the curvature equations \eqref{eq ut} satisfied by $u_s$ 
and $u_t$ ($0\leq s<t\leq 1$), we obtain
$$\Delta (u_t-u_s)= (K_s-K_t)e^{2u_t} +K_s(e^{2u_s}-e^{2u_t})$$
that we rewrite as~:
\begin{eqnarray}\label{eq ut-us}	
\Delta (u_t-u_s)=  (s-t)(1+K)e^{2u_t}+(-K_s)\,  (e^{2u_t}-e^{2u_s})\, .
\end{eqnarray}
Since the functions $u_t$ are uniformly bounded,
there exists a constant $m_0>0$, such that  one has 
$|u_t-u_s|\leq m_0\,|e^{2u_t}-e^{2u_s}|$  for all $s$, $t$ in $[0,1]$. 

The generalized maximum principle applied to $u_t-u_s$ combined with \eqref{eq ut-us}
ensures the existence of a constant $c$ such that $\| u_t-u_s\|_\infty\leq c\, |t-s|$ 
for every $s$, $t$ in $[0,1]$. 

Plugging this information into \eqref{eq ut-us} 
yields a similar bound for $\Delta (u_t-u_s)$,
and \eqref{eq u et grad u continu} follows  from the  Schauder estimates
\eqref{eq schauder 1} with $v=u_t-u_s$.\eproof\vspp

\noindent {\bf Remark } Since the curvature function $K$ is smooth, one could 
improve Lemma \ref{lem u borne et continu} and prove that all $u_t$ are smooth 
and that, for all $p\geq 2$ the maps $t\in [0,1] \to u_t\in C^p_{\rm loc}$ is continuous.
But the $p^{\rm th}$ derivatives of $u_t$ might not be bounded.

%62
\subsection{A continuous family of harmonic maps}

Recall that 
we have natural identifications $\partial_\infty S_t\simeq \SS^1$.
We fix an increasing quasi-symmetric homeomorphism $\varphi : \SS^1\to\SS^1$. 
In Chapter \ref{sec continuity}, we introduced the unique harmonic quasi-isometric map $h_t:\DD\to S_t$ with boundary value at infinity $\partial_\infty h_t=\varphi$. 

Here are the continuity properties of this  family of maps $h_t$ 
that we used in the proof of Theorem \ref{th main}. 

\bp
\label{prop continuity methodBplus} 
$(a)$ The map $t\in [0,1]\to h_t \in C(\DD,\DD)$ is continuous.\\ 
$(b)$ The map $t\in [0,1]\to J_t \in C_b(\DD,\RR)$ is continuous. 
\ep

\noi
This means that $\displaystyle\lim_{s\to t}d(h_s,h_t)=0$\; and\; 
$\displaystyle\lim_{s\to t}\|J_s-J_t\|_\infty=0$,\, for all $t\in [0,1]$.
\vspp 

\proof 
Assume this is not the case. Then there exist a sequence $(t_n)$ in $[0,1]$
and a sequence $(x_n)$ of points in $\DD$
such that 
\begin{eqnarray}
\label{eq htht jtjt}
\lim_{n\to\infty}d(h_{t_n}(x_n),h_t(x_n))>0
&{\rm or}&
\lim_{n\to\infty}|J_{t_n}(x_n)-J_t(x_n)|>0\, .
\end{eqnarray}

We want to get a contradiction by applying Lemmas \ref{lem seq surfaces}  and \ref{lem seq maps} 
to recentered surfaces and recentered harmonic maps. 
We thus choose  sequences $(\gamma_n)$ and $(\gamma'_n)$ of isometries of the hyperbolic disk
$\DD$ such that $\gamma_n(x_n)=0$ and $\gamma'_n(h_t(x_n))=0$. 
Let $S_n=(D,g_n)$ and $S'_n=(D,g'_n)$ be the conformal surfaces where $g_n=e^{2u_n}\gh$ and $g'_n=e^{2u'_n}\gh$ with
$$
u_n:= u_{t}\circ {\gamma'_n}^{-1}\;\; {\rm and}\;\;
u'_n:=u_{t_n}\circ {\gamma'_n}^{-1}\, .
$$
By Lemma \ref{lem seq surfaces}  we may assume, after extraction, that
the sequence $(u_n)$ converges 
to a $C^1$ function $u_\infty$ in the $C^1_{\rm loc}$ topology, 
and that the limit  $C^1$ metric space $S_\infty:=(D,e^{2u_\infty})$ is a
CAT space with pinched curvature $-b^2\leq K_\infty\leq -a^2<0$. 

By Lemma  \ref{lem u borne et continu}, one has 
$$
\displaystyle\lim_{n\to\infty}\|u'_n-u_n\|_\infty+\|\nabla u'_n-\nabla u_n\|_\infty=0\, .
$$ 
Hence the sequence $(u'_n)$  also converges  
in the $C^1_{\rm loc}$ topology to the function  $u_\infty$. 
We now introduce the sequence of maps 
\begin{eqnarray}
h_n:=\gamma'_n\circ h_{t}\circ \gamma_n^{-1} &:&\DD\to S_n\, ,\\
h'_n:=\gamma'_n\circ h_{t_n}\circ \gamma_n^{-1} &:&\DD\to S'_n\, .
\end{eqnarray}
These maps $h_n$ and $h'_n$ are harmonic and \eqref{eq htht jtjt} can be rewritten as
\begin{eqnarray}
\label{eq hnhn jnjn}
\lim_{n\to\infty}d(h_n(0),h'_n(0))>0
&{\rm or}&
\lim_{n\to\infty}|J_n(0)-J'_n(0)|>0\, ,
\end{eqnarray}
where $J_n$ is the Jacobian of $h_n$ and  $J'_n$   the Jacobian of $h'_n$. 
By Lemma \ref{lem hs ht}, all these maps $h_n$ and $h'_n$ are  uniformly Lipschitz and uniformly quasi-isometric.
Hence Lemma  \ref{lem seq maps} ensures that, after extraction,  
the sequences $(h_n)$ and $(h'_n)$ converge respectively, in the $C^2_{\rm loc}$ topology,  
to  harmonic quasi-isometric maps $h_\infty, h'_\infty:\DD\to S_\infty$. 

Since Lemma \ref{lem hs ht} also asserts that $d(h_n,h'_n)\leq 2 \, c_*$ for all $n$, the limit harmonic 
quasi-iso\-metric maps $h_\infty,h'_\infty:\DD\to S_\infty$ are within bounded distance from each other. 
Then 
the uniqueness theorem for quasi-isometric harmonic maps in \cite[\S 5]{BH18} ensures that $h_\infty=h'_\infty$. 
This contradicts \eqref{eq hnhn jnjn}.\eproof 
\vspp 

This also ends the proof of both Proposition \ref{prop continuity methodB} and  Theorem \ref{th main}.

%70
\section{ The injectivity theorem in constant curvature}
\label{sec litam}
\begin{quote}
This chapter is an appendix in which we prove the injectivity theorem \ref{th litam}  that we used as a starting point in the proof of our main theorem \ref{th main}.
\end{quote}

%71
\subsection{The Li-Tam-Markovic injectivity theorem}

\bt 
\label{th litam} 
Let $\DD$ be the hyperbolic disk.
Any harmonic  quasi-iso\-metric map $h:\DD\to\DD$  is a quasi-conformal harmonic diffeomorphism.
\et

\noi
This theorem is an output of Markovic solution of the Schoen conjecture in \cite{Marko2}.
It relies on a previous injectivity result of Li-Tam in \cite{LT}
when the boundary map of $h$ is smooth, which is Proposition \ref{pro fkdense} below. The proof of Li-Tam itself relies on the
Schoen-Yau injectivity theorem in \cite{SchoenYau78}.
 
We would like to give in this appendix a short new proof of Theorem \ref{th litam} 
that does not rely on this Schoen-Yau theorem and that uses instead
a continuity method combined with a simple topological fact (Lemma \ref{lem gorbit}). 
\vspp 

\proof
The proof will last till  the end of this appendix.
We know (see Section \ref{sec quasi-iso}) that the boundary value 
$\varphi=\partial_\infty h:\SS^1\to\SS^1$ is a $k$-quasi-symmetric homeomorphism of 
$\SS^1=\partial_\infty \DD$, where $k$ depends only on the constant $c$ 
of quasi-isometry of $h$.
For $k\geq 1$, we introduce the set 
\begin{eqnarray*}
\mk &=&\{ \hbox{ $k$-quasi-symmetric homeomorphism $\varphi:\SS^1\to\SS^1$}\}
\end{eqnarray*}
equipped  with the uniform distance 
$\displaystyle d(\varphi_1,\varphi_2)=\sup_{\xi\in\SS^1} |\varphi_1(\xi)-\varphi_2(\xi)|$.

We also know that, for all $\varphi$ in $\mk$, there exists a unique 
harmonic quasi-isometric map $h_\varphi:\DD\to\DD$ whose boundary map is $\varphi$. 
We want to prove that all these maps $h_\varphi$ are quasiconformal diffeomorphisms. 
This will follow from the next  Lemma \ref{lem diffeodense}, Proposition
\ref{pro fkclosed} and Proposition \ref{pro fkdense}.
\eproof

\bl
\label{lem diffeodense}
The $k$-quasi-symmetric $C^1$ diffeomorphisms are dense in $\mk$.
\el

\proof Choose a smooth approximation of unity $(\alpha_n)$ on $\SS^1$. 
For $\varphi$ in $\mk$, each function $\alpha_n*\varphi$ is a $k$-quasi-symmetric $C^1$
diffeomorphism
while the sequence $(\alpha_n*\varphi)$ 
converges uniformly to $\varphi$. 
\eproof

\bp
\label{pro fkclosed} Let  $\fk$ be the set of those $\varphi\in\mk$  such that  $h_\varphi$ is a quasi-conformal diffeomorphism. Then $\fk$ is a closed subset of $\mk$.
\ep

The proof of Proposition \ref{pro fkclosed} will be 
given in Section \ref{sec fkclosed}. It relies on continuity
properties of the boundary map $h\mapsto \partial_\infty h$
proven in Section \ref{sec conbou}.

\bp
\label{pro fkdense}
When $\varphi$ is a $C^1$ diffeomorphism of $\SS^1$, its quasi-isometric harmonic extension
$h_\varphi:\DD\to \DD$ is a quasi-conformal diffeomorphism.
\ep

The proof of Proposition \ref{pro fkdense} will be given in Section \ref{sec fkdense}. It uses a deformation $\varphi_t$ of $\varphi$ starting with the identity. Let $G$ be the group of isometries of $\DD$ acting on $\SS^1$.
The proof relies on the fact
that the only homeomorphisms which are limits of elements of   $G\varphi_t G$ belong to $G$.
This is Lemma \ref{lem gorbit} which will be proven in Section \ref{sec gorbit}.

%72
\subsection{Continuity of the boundary map}
\label{sec conbou}
Let $c>1$. Endow the space ${\cal Q}_c$ of $c$-quasi-isometric maps $f:\DD\to\DD$ with the topology of uniform convergence on compact sets, and the space ${\cal C}$ of continuous maps $\varphi:\SS^1\to\SS^1$ with the topology of  uniform convergence.

\bl
\label{lem conbou} 
The map $f\in {\cal Q}_c\to\partial_\infty f\in {\cal C}$ is continuous.
\el

\proof We use the quasi-invariance of the Gromov product under quasi-isometric maps.
We fix a point $0$ in $\DD$. For $n\in\NN\cup\{\infty\}$,
let $f_n\in {\cal Q}_c$ be $c$-quasi-isometric maps, with boundary values at infinity $\varphi_n$ . Assume that the sequence $(f_n)$ converges uniformly to $f_\infty$ on compact sets.
In particular, the quantity $R:=\sup_n d(f_n(0),0)$ is finite.
We want to prove that the sequence $(\varphi_n)$ converges uniformly to the boundary map $\varphi_\infty$ of $f_\infty$.

For $\xi\in\SS^1$, denote by $t\in\plu\to x_\xi^t\in\DD$ the geodesic ray with origin $0$ and endpoint $\xi$. By \cite[Proposition 5.15]{GhysHarp90}, there exists a constant $\lambda>1$ such that  the following lower bound for the Gromov product seen from  $0$
$$
(f_n(x_\xi^t),f_n(x_\xi^s))_{0} \geq (x_\xi^t, x_\xi^s)_0/\lambda-\lambda =t/\lambda -\lambda
$$
holds when $s\geq t>0$ and $n\in\NN\cup\{\infty\}$. Letting $s\to\infty$, we obtain
$$
(f_n(x_\xi^t),\varphi_n(\xi))_{0} \geq t/\lambda -\lambda\, 
$$
for $n\in\NN\cup\{\infty\}$. 
Since $\DD$ is $\delta$-hyperbolic for a constant $\delta>0$, each
Gromov product $(\varphi_n(\xi) ,\varphi_\infty(\xi) )_{0}$ is bounded below by
\begin{equation*}	
\min[ \, ( \varphi_n(\xi),f_n (x^t_\xi))_{0} \, ,\, 
(f_n (x^t_\xi), f_\infty (x^t_\xi)) _0 \, ,\,    (f_\infty (x^t_\xi), \varphi_\infty(\xi))_{0}\,  ] -2\delta
\end{equation*}
for every $\xi\in\SS^1$ and $n\in\NN$ (see \cite[Chap. 2]{GhysHarp90}).
The sequence $(f_n)$ converging uniformly to $f_\infty$ on compact sets
there exists, for all $t>0$, an integer $n_t\geq 1$ 
such that one has, for $n\geq n_t$ and $\xi\in\SS^1$,
$$d(f_n(x_\xi^t),  f_\infty(x_\xi^t))\leq 1
\; ,\;\;{\rm and \; hence}$$
$$
(\varphi_n(\xi) ,\varphi_\infty(\xi) )_{0}
\;\geq \;
\min  [t/\lambda -\lambda\, ; \, t/c-c-R-1/2] -2\delta\, .
$$
This proves that the sequence $(\varphi_n)$ converges uniformly to $\varphi_\infty$.
\eproof

%73
\subsection{A continuous inverse to the boundary map}
\label{sec fkclosed}

The following lemma is a variation of Lemma \ref{lem hs ht}. Fix $k\geq 1$.

\bl
\label{lem estimes hphi}
There exist a compact subset $L_k\subset\DD$ and  a constants $c_k$
such that
the harmonic quasi-isometric extension $h_\varphi$ of any  $\varphi\in\mk$
is $c_k$-quasi-isometric, the point  $h_\varphi (0)$ is in $L_k$, and the map
$h_\varphi$ is $c_k$-Lipschitz. 
\el

\proof 
We introduce  the Douady-Earle extension $f_\varphi:\DD\to\DD$ of $\varphi$ and
we recall some of their properties that can be found in J. Hubbard's book \cite[\S 5.1]{Hubb}. 
By definition, the image $f_\varphi(z)$ of $z\in \DD$ is the barycenter of the measure $\varphi_*(m_z)$
where $m_z$ is the visual measure on $\SS^1$ seen from $z$. 
This map $f_\varphi$ is smooth, and is $C_k$-quasi-isometric for some constant 
that depends only on $k$ (it is even $\delta_k$-quasi-conformal or some constant that depends only on $k$). 
The map $\varphi\to f_{\varphi}$ is continuous hence, since $\mk$ is compact,
the points $f_\varphi(0)$ belong to  a fixed compact set of $\DD$.

By the main result of \cite{Marko2} or \cite{BH15}, the distance
$d(h_\varphi, f_\varphi)$
is bounded by a constant $M_k$ that depends only on $C_k$. The first two claims follow.
The Lipschitz continuity of $h_\varphi$ then follows from the Cheng lemma \ref{lem cheng}.
\eproof

\bc
\label{cor phi uphi ctnu}
The map $\varphi\in\mk\to h_\varphi\in C^2(\DD,\DD)$ is continuous in the $C^2_{\rm loc}$ topology.
\ec

\proof Let $(\varphi_n)$ be a sequence in $\mk$ converging to $\varphi$. 
By Lemma \ref{lem estimes hphi}, the harmonic maps $h_n:=h_{\varphi_n}$ are 
uniformly locally bounded and uniformly Lipschitz. 
By Lemma \ref{lem seq maps}, after extraction, the sequence $(h_n)$ converges 
in the  $C^2_{\rm loc}$ topology to a harmonic quasi-isometric  map $h_\infty:\DD\to\DD$. 
To reach the conclusion, we need to prove that  such a limit $h_\infty$ is always equal to $h_\varphi$. 
Since the maps $h_n$ are uniformly quasi-isometric,  the continuity lemma \ref{lem conbou}
yields that
the limit $\varphi$ of the boundary maps $\varphi_n$ of $h_n$
must be the boundary map of $h_\infty$. This proves that  $h_\infty=h_\varphi$.
\eproof

\vsppp\noi{\bf Proof of Proposition \ref{pro fkclosed}}
Let $(\varphi_n)$ be a sequence in $\mk$ converging to $\varphi$ such that 
all the harmonic quasi-isometric extensions
$h_{\varphi_n}$ are quasiconformal diffeomorphisms. 
We want to prove that the harmonic map
$h_{\varphi}$ is also a quasiconformal diffeomorphism. 

Corollary \ref{cor phi uphi ctnu} ensures that  
the sequence $(h_{\varphi_n})$ 
converges to $h_\varphi$  in the $C^2_{\rm loc}$ topology. 
Lemma \ref{lem estimes hphi} ensures that  these maps $h_{\varphi_n}$ 
are uniformly Lipschitz.
Hence, by Proposition \ref{prop delta de m}, there exists
a uniform lower bound $j_*>0$ for the Jacobians  of all these 
harmonic quasi-isometric diffeomorphisms $h_{\varphi_n}$.
Therefore $h_\varphi$ is also a Lipschitz harmonic map whose Jacobian is bounded below by $j_*$.
Hence, by the injectivity criterion in Lemma \ref{lem continuity methodB}, 
the harmonic map $h_\varphi$ is also a quasiconformal diffeomorphism.
\eproof

%74
\subsection{Orbit closure in the group of 
homeomorphisms of $\SS^1$}
\label{sec gorbit}

Recall that $\DD$ is the hyperbolic disk and $\SS^1$ is its boundary at infinity.
Let $G$ be the group of isometries of $\DD$ acting on $\SS^1$. It is 
isomorphic to ${\rm PGL}(2,\RR)$. 

In order to prove Proposition \ref{pro fkdense} in the next section 
we will need the following lemma.

\bl 
\label{lem gorbit}
Let $\varphi_n$ be a sequence of $C^1$ diffeomorphisms of $\SS^1$ 
converging in the $C^1$\! topology to a $C^1$\! diffeomorphism $\varphi_\infty$ of $\SS^1$.\!
Let $\gamma_n$\! and $\gamma'_n$
be two unbounded sequences in $G$ such that the sequence
$\psi_n:=\gamma'_n\circ \varphi_n\circ\gamma_n^{-1}$ con\-verges 
to an homeomorphism $\psi_\infty$ of $\SS^1$. 
Then this limit $\psi_\infty$ belongs to $G$. 
\el

\proof 
We recall the Cartan decomposition $G=KA^+K$ of $G$ where 
$K$ is the group $PO(2,\RR)$ and 
$A^+=\{{\rm diag}(s,s^{-1}) \; {\rm with}\; s\geq 1\}$.
Since $K$ is compact,
we can assume that both $\gamma_n$ and $\gamma'_n$ are in $A^+$. 
We write 
$$
\gamma_n={\rm diag}(s_n^{1/2},s_n^{-1/2})
\;\;{\rm and}\;\;
\gamma'_n={\rm diag}({s'_n}^{1/2},{s'_n}^{-1/2})
$$
with both $s_n$ and $s'_n$ converging to $\infty$.
Here it will be convenient to use the identification $\SS^1\simeq \RR\cup\{\infty\}$ 
given by the upper half-plane model of $\DD$, so that, for $x$ in $\RR$, one has 
$\gamma_n(x)=s_nx$ and $\gamma'_n(x)=s'_nx$.

We  notice that $\varphi_\infty(0)=0$.
Indeed if this were not the case, we would have $\psi_\infty(x)=\infty$ for all $x\in \RR$,
contradicting the injectivity of $\psi_\infty$.

Similarly we have $\psi_\infty(\infty)=\infty$.
Indeed if this were not the case, we would have $\varphi_\infty(x)=0$ for all $x\in \RR$,
contradicting the injectivity of $\varphi_\infty$.

Since the sequence $\varphi_n$ converges in the $C^1$ topology to $\varphi_\infty$,
we can write for all $n\geq 1$ and all $x\in \RR$ with $|x|\leq 1$
\begin{equation}
\varphi_n(x)=\alpha_n +(\beta_n +r_n(x))x
\;\;\;
\mbox{with}\;\;\;
\lim_{x\to 0}\sup_{n\in \NN}|r_n(x)|=0\, .
\end{equation}
Since $\varphi_\infty(0)=0$ and $\beta_\infty:=\varphi'_\infty(0)$ is non zero, one has 
\begin{equation}
\lim_{n\to\infty} \alpha_n=0
\;\;{\rm and}\;\;
\lim_{n\to\infty}\beta_n=\beta_\infty>0\, .
\end{equation}
Therefore we can write for all $n\geq 1$ and all $x\in \RR$ with $|x|\leq s_n$
\begin{equation}
\psi_n(x)=s'_n\alpha_n +(\beta_n +r_n(\tfrac{x}{s_n}))\tfrac{s'_n}{s_n}x
\;\;\;
\mbox{with}\;\;\;
\lim_{n\to \infty}|r_n(\tfrac{x}{s_n})|=0\, .
\end{equation}

Since the sequences $\psi_n(0)$ and $\psi_n(1)$ converge, the following  limits  exist 
\begin{equation}
\alpha'_\infty:=\lim_{n\to\infty} s'_n\alpha_n\in \RR
\;\;{\rm and}\;\;
\beta'_\infty:=\lim_{n\to\infty}\beta_n\tfrac{s'_n}{s_n}>0\, ,
\end{equation}
Hence one has $\psi_\infty(x)=\alpha'_\infty+\beta'_\infty x$ for all $x\in \RR$,
and $\psi_\infty$ belongs to $G$. \eproof
\vspp  

\noi
{\bf Remark} - As can be seen in the proof, 
the assumption on $\psi_n$ can be weakened: it is sufficient to assume that
there are three points $\xi_0$, $\xi_1$, $\xi_\infty$ in $\SS^1$ whose images
$\psi_n(\xi_0)$, $\psi_n(\xi_1)$, $\psi_n(\xi_\infty)$ converge to three distinct points.
This ensures that the sequence $\psi_n$ converges uniformly
to an element $\psi_\infty$ of $G$.\\
- However, it is important to assume
that the limit $\varphi_\infty$ is of class $C^1$ 
and that the convergence $\varphi_n\to\varphi_\infty$
is in the $C^1$ topology.
\vspp  

Here is a direct corollary of Lemma \ref{lem gorbit} 
in the spirit of  \cite{BHAutosimilar}.

\bc
\label{cor gorbit}
For all $C^1$ diffeomorphism $\varphi$ of $\SS^1$, one has the equality\\ 
\centerline{$\overline{G\varphi G}\,\cap\, {\mathcal H}omeo(\SS^1)
\;=\; 
G\varphi G \,\cup\, G$.}
\ec

%75
\subsection{When the boundary map is a $C^1$ diffeomorphism}
\label{sec fkdense}

We now conclude the proof of Theorem \ref{th litam} by giving the last argument:
\vspp 

\noi
{\bf Proof of Proposition \ref{pro fkdense} } Let $\varphi$ be a 
$C^1$ diffeomorphism of $\SS^1$.
We want to prove that the harmonic  quasi-isometric extension $h_\varphi$
of  $\varphi$ is a quasi-conformal diffeomorphism.
For convenience we identify here $\SS^1$ with $\RR/2\pi\ZZ$. 
For $t\in[0,1]$, we introduce the $C^1$ diffeomorphism $\varphi_t$ 
given by 
$$
\varphi_t(\xi)=\xi+(\varphi(\xi)-\xi)\, t
\;\;\;
\mbox{ for all $\xi$ in $\SS^1$.}
$$
This is well defined since the map $\xi\to \varphi(\xi)-\xi$ lifts as a map 
from $\SS^1$ to $\RR$. 

We argue as in Section \ref{sec strategy}. 
For $t\in [0,1]$ we introduce the  harmonic quasi-isometric extension $h_t=h_{\varphi_t}:\DD\to\DD$
of $\varphi_t$. 
Let $A$ be the set of  parameters $t\in [0,1]$ for which
$h_t$ is a quasi-conformal diffeomorphism. 
By the injectivity criterion of Lemma \ref{lem continuity methodB}, one has
$$
A=\{t\in [0,1]\mid \inf_{z\in \DD}J_t(z)>0\}\,
$$
where $J_t$ is the Jacobian of $h_t$.
We want to prove that $1\in A$.
We already know that $0\in A$  because $h_0$ is the identity.
Since the maps $\varphi_t$ are 
uniformly quasi-symmetric, 
Proposition \ref{pro fkclosed} tells us  that $A$ is closed.
Therefore it is enough to check that $A$ is open.

Assume by contradiction that there exists a sequence $t_n\not\in A$ converging to 
$t_\infty \in A$.  
By assumption there exists a sequence $(z_n)$ in $\DD$ such that 
$\displaystyle\liminf_{n\to\infty} J_{t_n}(z_n)\leq 0$.
After extraction we are in one of the two cases:
\vspp

{\bf First case } The sequence $(z_n)$ converges to a point $z_\infty\in \DD$.\\
Since the maps $\varphi_t$ are 
uniformly quasi-symmetric,
Corollary \ref{cor phi uphi ctnu} ensures that
the map 
$t\in [0,1]\to h_t\in C^2(\DD,\DD)$ is continuous in the $C^2_{\rm loc}$ topology.
Therefore, one has 
$
\displaystyle
J_{t_\infty}(z_\infty)=
\lim_{n\to\infty} J_{t_n}(z_n)\leq 0
$,
and $t_\infty$ is not in $A$. Contradiction.
\vspp  

{\bf Second case } The sequence $(z_n)$ goes to infinity.\\
To simplify
we set $\varphi_n=\varphi_{t_n}$ and $h_n=h_{t_n}$ for  all $n\in \NN\cup \{\infty\}$.
By Lemma \ref{lem estimes hphi}, the sequence $h_n(z_n)$  goes to infinity. 
We choose  sequences $(\gamma_n)$ and $(\gamma'_n)$ in $G$ 
with $\gamma_n(z_n)=0$ and $\gamma'_n(h_n(z_n))=0$. 
We  introduce the  harmonic maps 
\begin{eqnarray*}
h'_n:=\gamma'_n\circ h_n\circ \gamma_n^{-1} &:&\DD\to \DD\,  
\end{eqnarray*}
and their boundary values $\psi_n:=\gamma'_n\circ \varphi_n\circ\gamma_n^{-1}$.
By construction, one has 
\begin{equation}
\label{eqn hn0jn0}
h'_n(0)=0
\;\;{\rm and}\;\; 
\liminf_{n\to\infty}J'_n(0)\leq 0\, ,
\end{equation}
where $J'_n$ is the Jacobian of $h'_n$.
Moreover by Lemma \ref{lem estimes hphi},
these maps $h'_n$ are uniformly Lipschitz. Therefore, after extraction, 
they converge in the $C^2_{\rm loc}$ topology to a harmonic quasi-isometric map $h'_\infty$.
By the continuity lemma \ref{lem conbou}, the sequence of boundary maps $\psi_n$
converge to the boundary map $\psi_\infty$ of $h'_\infty$.
Now, by Lemma \ref{lem gorbit}, this limit $\psi_\infty$ belongs to $G$. 
Therefore the harmonic map  $h'_\infty$ is an isometry  and its Jacobian is 
$J'_\infty\equiv 1$. This contradicts \eqref{eqn hn0jn0}.
\eproof

%90
%\small{	\bibliographystyle{plain}	\bibliography{harmonic4}}
\small{

}
\vspace{2em}

{\small\noindent Y. Benoist  \& D. Hulin,\; 
	CNRS \& Université Paris-Saclay \\ Laboratoire de mathématiques d’Orsay, 91405, Orsay, France
	
	\vspace{0.3em}
\noindent yves.benoist@u-psud.fr \;\&\; dominique.hulin@u-psud.fr}

\end{document}